\documentclass[a4paper,10pt]{article}
\def \Z {{\mathbf {Z}}}

\def \P {{\mathbf {P}}}

\def \m {{ \mu}}

\def\eps{\varepsilon}
\title{ Ограниченные эргодические конструкции,\\ дизъюнктность и слабые пределы степеней }
\author{В.В. Рыжиков\footnote{\large Работа поддержана грантом  НШ-5998.2012.1.
}}
\date{ 28.03.13}
\usepackage[T1]{fontenc} 
\usepackage[cp1251]{inputenc} 
\usepackage[russian]{babel}

\begin{document}
\Large
\maketitle

\begin{abstract}{\large Статья посвящена свойству   дизъюнктности степеней вполне эргодической ограниченной конструкции ранга  один и  некоторым обобщениям этого результата. Рассмотрены  приложения к задаче  независимости функции Мёбиуса от последовательности, индуцированной ограниченй конструкцией.

Ключевые слова:  Эргодическая степень  преобразования, конструкция ранга один,  дизъюнктность динамических систем, функция Мёбиуса.  	}
\end{abstract}

Интерес к озаглавленой тематике связан со  следующим наблюдением.  Ограниченные конструкции ранга один при условии, что все их ненулевые степени эргодичны,  обладают нетривиальными  слабыми  пределами степеней.  Это влечет за собой дизъюнктность (в смысле работы \cite{F})  степеней конструкций,  что в силу результатов \cite{BSZ} приводит к  независимости  ограниченных конструкций  от  функции Мёбиуса. 
Таким образом, задача о дизъюнктности степеней преобразований, ранее рассматриваемая 
как специальная задача в рамках теории самоприсоединений,
получила интересное приложение. 

 Известная гипотеза  \cite{S} утверждает для  строго эргодического  гомеоморфизма $S:X\to X$ с нулевой топологической энтропией выполнение свойства 
 $$\sum_{i=1}^N f(S^{i} x)\m(i)=o(N),$$
где $f\in C(X)$,  а   $\m$ -- функция Мёбиуса.   Это свойство в статье называется  независимостью от функции Мёбиуса (в цитированной литературе оно называется свойством дизъюнктности или ортогональности).

 В  \cite{B} Бурген исследовал  произведения Рисса   и    подтвердил  упомянутую  гипотезу Сарнака для ограниченных конструкций ранга один.  
В \cite{R} методом слабых пределов степеней  показано, что для вполне эргодической  ограниченной конструкции ее положительные  степени дизъюнктны.   Это упростило доказательство упомянутого результата Бургена. 
В работе \cite{ALR}    авторы доказали спектральную дизъюнктность степеней
слабо перемешивающей  конструкции ранга один,  обладающей  
так называемой  неплоской  ограниченной рекуррентностью.
  Отметим, что спектральную дизъюнктность степеней преобразования иногда можно устанавливать методом слабых пределов
даже при отсутствии нетривиальных пределов (см. \cite{R3}).  Для наших целей достаточно установить только дизъюнктность
степеней.  Отметим, что недавно А.А. Приходько построил поток ранга 1 с  лебеговским спектром \cite{P}.  Преобразования, входящие в такой  поток, обладают спектрально изоморфными степенями. Однако,  их степени дизъюнктны в смысле работы \cite{F}, что было замечено в   \cite{R1}.  

Цель этой заметки  -- изложить, следуя
идеям \cite{R}, \cite{ALR},  простое  доказательство 
дизъюнктности степеней упомянутых конструкций.
Также будет показано, как редуцировать эргодический случай 
к вполне эргодическому, чтобы  установить  свойство независимости ограниченых конструкций ранга один от функции Мёбиуса.

\section{Достаточное условие дизъюнктности степеней преобразования}
Обратимые  сохраняющие меру преобразования вероятностного пространства Лебега $(X,\nu)$ в статье для краткости называются \it преобразованиями. \rm
 Ниже мы пользуемся  терминологией  марковских (стохастических)  сплетающих операторов (см., например,   \cite{R1}).  Преобразования  $S$ и $T$ 
(как операторы в  пространстве  $L_2(X,\nu), \ \nu(X)=1$) называются \it дизъюнктными, \rm если для любого марковского оператора $J$ равенство 
 $SJ=JT$  влечет за собой  $J=\Theta$, где $\Theta$  -- ортопроекция  на пространство констант  в  $L_2(X,\nu)$. Это определение дизъюнктности эквивалентно класcическому определению $\cite{F}$, использующему понятие присоединения (joining).

Преобразование  $T$ называется \it вполне эргодическим, \rm если все его ненулевые степени эргодичны. Это эквивалентно тому, что для любого марковского оператора $J$ и любого  $i>0$ равенство $T^iJ=J$  влечет за собой   $J=\Theta$.

Марковский оператор  $J$, удовлетворяющий  условию сплетения  $SJ=JT$,   будем  называть  \it неразложимым, \rm если  его нельзя представить в виде выпуклой суммы различных марковских операторов, сплетающих $S$ и $T$.   
\vspace{5mm}

\bf Лемма  1. \it  Пусть  $S,T$  -- вполне эргодические пребразования, а   $J\neq \Theta$ -- неразложимый оператор,
удовлетворяющий условию сплетения  $$S^qJ=JT^p,$$
где $q,p$ взаимно просты.
Если выполнено   $$Q(S)J=JP(T),$$
 где  $ Q(S)=\sum_{i} a_iS^i,  \ \ P(T)=\sum_{j} b_j T^j,  {}\hspace{5mm}{} \, \, 
  \ \
\sum_{i} a_i =1=\sum_{i} b_j,     \  a_i,b_j\geq 0,$\\
то  ряды $Q$  и  $P$ являются  $\frac p q$ -подобными:  найдется ряд $R$ такой, что
$$ Q( S)=R(S^q) ,  {}\hspace{5mm}{} P(T)=R(T^p) . $$

\rm 
Доказательство.  Заметим, что  операторы вида  $S^kJ, JT^m$  неразложимы, так  $S^k, T^m$ обратимы, а $J$ неразложим. Тогда из равенства $Q(S)J=JP(T)$  вытекает, что  члены  ряда
$\sum_{i} a_iS^iJ$ совпадают с членами ряда  $\sum_{j} b_j JT^j$. Следовательно,   при     $a_i\neq 0$ для некоторого $j$  выполнено  $b_j\neq 0$  и 
 $$S^iJ=JT^j.$$
Если  $i\neq 0$,  то  $j\neq 0$, так как из эргодичности $S^i$ при условии  
$S^iJ=J$  получим
$J=\Theta.$
Имеем  $S^{i}J=JT^{j},\ \ ij\neq 0.$
Тогда $S^{iq}J=JT^{jq},$  а
из $S^{q}J=JT^{p}\ $ следует
  $S^{iq}J=JT^{ip}.$
Таким образом, 
$J=JT^{jq-ip}.$   Но  преобразование $T$ вполне эргодично,
значит,  
$jq=ip.$
Так как   $q,p$ взаимно просты,  получаем 
$i=qr,  \ \    j=pr.$
Тем самым мы показали, что  
\begin{center} 
$a_i=b_{iq/p}$, \ $\sum_{r} a_{qr} =1$, \ 
$\sum_{r} b_{pr} =1$.
\end{center}
 Это означает то, что ряды    $Q$ и $P$ являются  $\frac p q$ -подобными.

Из леммы 1 непосредственно вытекает  

\vspace{2mm}
\bf  Лемма 2. \rm  \it Пусть вполне эргодическое преобразование $T$ обладет слабыми пределами вида 
$T^{qn_i}\to Q(T),    \ \   T^{pn_i}\to P(T).$  Если ряды  $Q(T)$ и  $P(T)$ не являются   $\frac p q$-подобными,  то 
 $T^q$ и $T^p$  дизъюнктны. \rm

\section{ Дизъюнктность степеней частично ограниченных конструкций}

Конструкция ранга 1 задается  следующим  параметрами: числом  $h_1$ (высота начальной башни), последовательностью  $r_j$ (число колонн, на которое разрезается башня на этапе с номером $j$)  и  последовательностью массивов из высот надстроек над колоннами 
$ (s_j(1), s_j(2),\dots, s_j(r_j-1),s_j(r_j)).$ 
Высота башни с номером  $j$ вычисляется по формуле
$h_{j+1}+1=(h_j+1)r_j +\sum_{i=1}^{r_j}s_j(i).$
Подробное описание конструкций ранга 1 можно найти в цитируемых работах. 
 Если все параметры  $r_j$  и 
$ s_j(i)$  ограничены,  то конструкция называется \it ограниченной. \rm
\vspace{2mm}

\bf Теорема  1.  \it Пусть $T$  --  вполне эргодическая  ограниченная конструкция.  Тогда степени  $T^p$ и $T^q$  дизъюнктны для всех положительных  чисел $q\neq p$. \rm

\vspace{2mm}
Несложное доказательство этой теоремы приведено в \cite{R} (см.  \S 7).
В \cite{ALR} авторы рассмотрели следующее интересное обобщение
ограниченных конструкций.  
Назовем  конструкцию \it частично ограниченной, \rm если  параметры конструкции  ограничены   для этапов  $j$   из некоторого подмножества 
натурального ряда, содержащего сколь угодно большие целочисленные  интервалы.  Из этой системы интервалов можно выбрать такую подсистему интервалов $J_k=\{j_k, j_k+1, \dots , j_k+l_k\}$, что  параметры  конструкции
совпадают  для всех  этапов вида  $ j_{k}+m$  при $k\geq m$.
 В работе \cite{ALR} эквивалентное  свойство называлось \it рекуррентной ограниченностью.  \rm  Если же величина высот надстроек 
стабилизируется на некотором постоянном значении, т.е.
$s_j(i)= s$  для всех $j\in J_k$  и $i\leq r_j$, то такие конструкции называют \it рекуррентно плоскими. \rm

Рассмотрим 
взаимосвязи между слабыми пределами степеней
рекуррентно ограниченных преобразований.  
Далее  используем следующие обозначения   $$H_j=-h_j - s_j,  \ \ s_j= min \{s_j(1),s_j(2), \dots, s_j(r_j-1)\}$$ и 
$$ P_{d,m} (T):= \lim_{k\to\infty}   T^{dH_{j_k+m}}.$$
 
Стандартные  вычисления слабых пределов степеней  преобразований ранга 1  приводят к следующим соотношениям:
$$ P_{d,m} (T) =  \dots +  a_{d,m}T^{N(d,m)}P_{1,m+1}, \eqno (\ast)$$
$$ P_{1,m} (T) = aI+  \dots, \ \ \ a>0.$$

Покажем, как из этих соотношений вытекает неограниченность    или же плоское поведение  надстроек.
Для  $P(T)=\sum_z a_zT^z$  определим  $\P\subset \Z$,
положив  $$ \P=\{z\in \Z :   a_z>0\}.$$
Из леммы 2  имеем
$$ \P_{p,m}\subset   p\Z.$$
Отметим, что из  $ P_{1,m} (T) = aI+  \dots $  (тем самым $0\in  \P_{1,m+1}$) в силу $(\ast)$ вытекает, что  
$$ \P_{1,m+1}\subset   p\Z.$$

В силу того, что любая  сумма чисел $s_j(i)$ делится на  $p$,  при условии, что  каждое из
чисел $s_j(i)$ делится на  $p$,  получим  
\\
$  \P_{q,m+1}\subset ( \P_{1,m+1}+ \P_{1,m+1}+\dots +\P_{1,m+1})$\ \ ($q$ \ слагаемых). 
Это влечет за собой 
$$ \P_{q,m+1}\subset   p\Z.$$
Из  леммы 2  следует
$$ \P_{p,m+1}=\frac p q  \P_{q,m+1}\subset   p^2\Z.$$
Итак
$$ \P_{1,m+2}\subset   p^2\Z,$$
 повторяя рассуждение, получим 
$ \P_{1,m+3}\subset   p^3\Z,  \ \ \P_{1,m+4}\subset   p^4\Z, \dots .$

Тем самым для частично ограниченных конструкций  мы показали, что для любого натурального  $m$ при  достаточно большом  $k$  ограниченные разности   $s_{j_k +m}(i)- s_{j_k +m}('i)$ делятся на $p^m$.  Условие ограниченности
параметров конструкций запрещает это для больших значений $p^m$.

Получили для всех  $m$ таких, что $p^m>s$ плоское поведение надстроек: 
$$ s_{j_k +m}(1)= s_{j_k +m}(2)=\dots= s_{j_k +m}(r_{j_k +m}-1):=s_{k,m}.$$
Нетрудно видеть, что для фиксированных $m, m'$ локальные константы склеиваются: $ s_{k,m}=s_{k,m'}$  для всех больших $k$. Действительно,  если константы различны, то  в слабом замыкании возникает нетривиальный ограниченный полином,  который 
в силу сказанного выше  запрещен.   Подведем итог.
\vspace{3mm}

\bf Теорема  2.  \it Вполне эргодическая частично ограниченная конструкция, у которой степени  $T^p$ и $T^q$ не дизъюнктны для некоторых взаимно простых чисел $q, p$ является рекурентно плоской  конструкцией.\rm

Эта теорема вытекает из результатов \cite{ALR}. Мы привели 
упрощенное доказательство.

 В случае  плоских  конструкций возникает дополнительная возможность: после серии плоских надстроек разрешаются  неограниченные  надстройки.   Для широкого класса таких надстроек из условия сплетения $T^qJ=JT^p$ можно получить соотношение вида  
$$((1-q\eps) I + q\eps\Theta)J=J((1-p\eps) I + p\eps \Theta),
\eqno {\bf (I,\Theta)}$$
где $\eps>0$  -- некоторое маленькое число. Но это  при $q\neq p$ очевидным образом приводит к дизъюнктности степеней $T^q$ и $T^p$,  т. е. к равенству $J=\Theta$.

Имееются  слабо перемешивающие преобразования ранга 1  с изоморфными степенями \cite{A},\cite{D},\cite{R2}.
Такие преобразования    не могут обладать слабыми пределами степеней, которые фигурируют в формулах  ${\bf (I,\Theta)}$ и  $\bf (I,P)$ (см. ниже).
\vspace{3mm}
 
 \vspace{5mm}

\Large
\section{ Независимость функции Мёбиуса и ограниченой конструкции}

В \cite{B}  Бургеном было доказано следующее утверждение.

 \vspace{3mm}
\bf Теорема 3. \it  Ограниченные конструкции ранга 1 независимы от функции Мёбиуса.\rm
\vspace{3mm}

Доказательство.  
Ограниченная  конструкция принадлежит одному из перечисленных классов: 

  1.  системы с дискретным рациональным спектром (одометры), 

  2.  плоские слабо перемешивающие конструкции, 

   3.  неплоские  слабо перемешивающие конструкции,

  4.  неплоские конструкции  с конечным компактным фактором (некоторая степень этой конструкции является прямой суммой слабо перемешивающих  ограниченных слабо перемешивающих конструкций).
 Таким образом, доказательство разбивается на 4 части.

1. Для одометров независимость от  функции Мёбиуса  вытекает из известного факта  
$  \sum_{i=1}^N  \m (pi) =o(N).$
Дело в том, что для любого $n$ одометр  является циклической перестановкой 
$p>n$  множеств $E, TE,\dots, T^{p-1}E$, индикаторы которых в виде линейной комбинации
близки к заданной непрерывной функции. Поэтому задача сводится к рассмотрению индикаторов $\chi_E$ этих множеств. 
А для них все известно:  для $x\in E$
$$  \sum_{i=1}^N  \chi_E (T^ix)\m (i) = \sum_{0<i\leq N/p}  \m (pi) =o(N). $$

2.  
  Как было показано в \cite{R}(\S 7),  если ограниченная конструкция не является одометром,  но обладает плоскими надстройками (является плоско рекуррентной),  то и в  этом случае возникают  нетривиальные полиномиальные пределы  вида $(1-m\eps) I + m\eps P$.  А именно,  из 
условия сплетения  $T^qJ=JT^p$  получим
$$((1-q\eps) I + q\eps P)J=J((1-p\eps) I + p\eps P)  \eqno {\bf (I,P)}$$
для некоторого полинома  $P\neq I$ и малого $\eps>0$, что приводит к $J=JT^k$, $k\neq 0$,  следовательно,  если $T^k$  -- эргодическая степень, то $J=\Theta$.  
Это влечет за собой (в силу \cite{BSZ}) независимость "процесса" \ $T$ от  функции Мёбиуса.

3.  В этом случае применяем  теорему 2.

4.   Осталось изучить случай, когда некоторая степень
преобразования $T$  не является эргодической.
Если ограниченная конструкция  $T$ не является одометром и не является вполне эргодической,  то $T$ как оператор обладает собственным значением $\lambda$ таким, что для некоторого (минимального) числа   $d>0$ выполнено
 $\lambda^d=1$,  причем $\lambda$  порождает группу всех собственных значений оператора $T$. 
Этот факт есть следствие наличия нетривиального полинома $P(T)=\sum_i a_i T^i \neq I$ как слабого предела степеней.
Действительно, из
$\left|\sum_i a_i \lambda^i\right|=1$ следует, что при  $a_k\neq 0$ выполнено
$\lambda^k=1$.  Положим $d=min\{k-k' :\ a_ka_{k'}\neq 0, k\neq k'\}$.
 
Такая группа собственных значений возникает только в том случае, если  автоморфизм  $T$ циклически переставляет некоторые множества 
$E, TE,\dots, T^{d-1}E$, причем ограничение   $S=T^d |E$ является слабо перемешивающим автоморфизмом (имеет непрерывный спектр) на пространстве $E$.

Число  $d$ обязано делить высоты колонн с надстроенными над ними этажами,
т. е.  $d$  делит числа $h_j +s_j(i)$  для всех достаточно больших $j$. Действительно,  пусть этаж башни на  "бесконечно большом" шаге $j$ с бесконечно малой относительной погрешностью содержится в множестве $E$.  Если бы нашлась такая колонна с номером $i$, что  $h_j +s_j(i)=dM+q$, $0<q<d$,  то получим, что часть множества  $E$ меры больше, чем $\nu(E)/r -\eps$ (число $r$ ограничивает число колонн) оказалась бы  вне множества $E$.  Но это противоречит  инвариантности $E$ относительно степени $T^d$.

Таким образом,  $S=T^d |E$ является ограниченной слабо перемешивающей конструкцией.  Но ранее мы уже доказали дизъюнктность степеней для нее.  Поэтому $S$ независима от функции Мёбиуса:  для  $f\in C(E)$ (мы рассматриваем  E как отрезок) и $n>0$  имеем
$$\sum_{i=1}^N f(S^{ni} x)\m(i)=o(N).$$
Рассмотрим  $f$  как функцию на  $X$, доопределив ее нулем вне множества $E$.
Пусть  $d$ простое. Для функции Мёбиуса  выполнено  $\m(d)=-1$, 
$\m(d^2k)=0$ и $\m(dk)= \m(d)\m(k)$ при $k\perp d$. Для $x\in E$ получаем
$$\sum_{i=1}^N f(T^{i} x)\m(i)=\sum_{0<k\leq N/d} f(T^{dk} x)\m(dk)=
\sum_{0<k\leq N/d} f(S^{k} x)\m(d)\m(k) - 
$$
$$-\sum_{0<m\leq N/d^2} f(S^{dm} x)\m(d)\m(dm)=
\frac {o(N)} {d}  + 
\sum_{0<m\leq N/d^2} f(S^{dm} x)\m(dm)= $$
$$= \frac {o(N)} {d}  + 
\sum_{0<m\leq N/d^2} f(S^{dm} x)\m(d)\m(m)+$$
$$+
\sum_{0<n\leq N/d^3} f(S^{d^2n} x)\m(dn) \leq \frac {o(N)} {d}  + \frac {o(N)} {d^2}  + \frac {N\|f\|} {d^3}.
$$
Аналогично,  для любого $M$ для всех достаточно больших $N$
получим 
$$\left |\sum_{i=1}^N f(T^{i} x)\m(i)\right|\leq \frac {N\|f\|} {d^{M}}; \ \ \ \sum_{i=1}^N f(T^{i} x)\m(i)=o(N).
$$

Теперь представим $F\in C(X)$  как $F=\sum_{i=0}^{d-1} f_i$, где
$supp f_i \subset T^iE$. Получаем
$$\sum_{i=1}^N F(T^{i} x)\m(i)=o(N).$$
Мы показали, что простые (prime) расширения  сохраняют свойство
независимости от функции Мёбиуса.
Если же  $d$ не является простым числом,  выполняем последовательно серию простых расширений в соответствии с разложением
числа $d$ на простые множители.

Из теоремы  2 и  описанной редукции  эргодического случая к вполне эргодическому получаем следующий факт \cite{ALR}:
\it  частично ограниченная конструкция  является плоско рекуррентной или обладает свойством независимости от функции Мёбиуса.\rm
\vspace{2mm}

Автор благодарит Ж.-П. Тувено (J.-P. Thouvenot) за стимулирующие беседы, без которых эта статья вряд ли 
была бы написана, и выражает признательность С.В. Тихонову за полезные замечания. 
\large

\vspace{5mm}
Москва, МГУ им. М. В. Ломоносова, \\
механико-математический факультет \\
E-mail: vryzh@mail.ru
\end{document}